\documentclass{article}
\usepackage[latin2]{inputenc}

\usepackage{amsmath}
\usepackage{amssymb}
\usepackage{graphicx}

\def\dom{{\rm dom}}
\def\proof{\noindent{\sc Proof. }}
\def\qed{\hspace{\stretch{1}}$\square$\medskip}

\DeclareMathOperator{\graph}{graph}

\DeclareMathOperator{\supp}{supp}

\DeclareMathOperator{\ran}{ran}

\newcommand{\ZZ}{\mathbb{Z}}

\newcommand{\si}{\sigma}

\newcommand{\de}{\delta}
\renewcommand{\epsilon}{\varepsilon}
\newcommand{\eps}{\varepsilon}

\newcommand{\sm}{\setminus}

\newcommand{\beeq}{\begin{equation}}
\newcommand{\eeeq}{\end{equation}}
\def\su{\subset}

\begin{document}

\title{Haar null sets without $G_\delta$ hulls}

\author{M\'arton Elekes\thanks{Partially supported by the
Hungarian Scientific Foundation grants no.~72655, 83726.}
\ and
Zolt\'an Vidny\'anszky\thanks{Partially supported by the
Hungarian Scientific Foundation grant no.~104178.}}

\insert\footins{\footnotesize{MSC codes: Primary 03E15, 54H05; Secondary 54H11, 28A99, 03E17, 22F99}}
\insert\footins{\footnotesize{Key Words: non-locally compact Polish group, separable Banach space, Haar null, Christensen, shy, prevalent, additivity, add, Mycielski, $G_\delta$, hull, coanalytic, universally measurable, Fremlin, Problem GP}}

\maketitle

\abstract{Let $G$ be an abelian Polish group, e.g.~a separable Banach space. A
subset $X \su G$ is called \emph{Haar null (in the sense of Christensen)} if
there exists a Borel set $B \supset X$ and a
Borel probability measure $\mu$ on $G$ such that $\mu(B+g)=0$ for every $g \in
G$. The term \emph{shy} is also commonly used for Haar null, and
co-Haar null sets are often called \emph{prevalent}.

Answering an old question of Mycielski we show that if $G$ is \emph{not}
locally compact then there exists a Borel Haar null set that is not contained
in any $G_\delta$ Haar null set. We also show that $G_\delta$ can be replaced
by any other class of the Borel hierarchy, which implies that the additivity
of the $\si$-ideal of Haar null sets is $\omega_1$. 

The definition of a \emph{generalised Haar null set} is obtained by replacing the
Borelness of $B$ in the above definition by universal measurability.
We give an example of a generalised Haar null set that is not Haar null, more precisely
we construct a coanalytic generalised Haar null set without a Borel Haar null hull.
This solves Problem GP from Fremlin's problem list. Actually,
all our results readily generalise to all Polish groups that admit a two-sided
invariant metric.}
 

\newtheorem{theor}{jejj}[section]
\newtheorem{theorem}[theor]{Theorem}
\newtheorem{co}[theor]{Corollary}
\newtheorem{all}[theor]{Proposition}
\newtheorem{lemma}[theor]{Lemma}
\newtheorem{krd}[theor]{Question}
\newtheorem{example}[theor]{Example}
\newtheorem{defin}[theor]{Definition}
\newtheorem{claim}[theor]{Claim}
\newtheorem{remark}[theor]{Remark}
\newtheorem{fact}[theor]{Remark}
\newcommand{\thm}[2]{\begin{theorem}\label{#1}#2\end{theorem}}
\newcommand{\cor}[2]{\begin{co}\label{#1}#2\end{co}}
\newcommand{\prop}[2]{\begin{all}\label{#1}#2\end{all}}
\newcommand{\lem}[2]{\begin{lemma}\label{#1}#2\end{lemma}}
\newcommand{\prob}[2]{\begin{pr}\label{#1}#2\end{pr}}
\newcommand{\ex}[2]{\begin{example}\label{#1}#2\end{example}}
\newcommand{\defi}[2]{\begin{defin}\label{#1}#2\end{defin}}
\newcommand{\rem}[2]{\begin{remark}\label{#1}#2\end{remark}}
\newcommand{\fct}[2]{\begin{fact}\label{#1}#2\end{fact}}
\newcommand{\que}[2]{\begin{krd}\label{#1}#2\end{krd}}
\newcommand{\clm}[2]{\begin{claim}\label{#1}#2\end{claim}}
\newcommand{\R}{\mathbb{R}}
\newcommand{\om}{\omega}
\newcommand{\md}{\models}
\newcommand{\bbb}{\mathbb}

\section{Introduction}

Throughout the paper, let $G$ be an abelian Polish group, that is, an
abelian topological group that is separable and admits a complete metric
(the group operation will be denoted by $+$ and the neutral element by $0$). It
is a well-known result of Birkhoff and Kakutani that any metrisable group admits
a left invariant metric \cite[1.1.1]{kechbeck}, which is clearly two-sided invariant for abelian
groups. Moreover, it is also well-known that a two-sided invariant metric on a
Polish group is complete \cite[1.2.2]{kechbeck}. Hence from now on let $d$ be a fixed complete
two-sided invariant metric on $G$. For the ease of notation we will restrict our attention to
abelian groups, but we remark that all our results easily
generalise to all Polish groups admitting a two-sided
invariant metric.

If $G$ is locally compact than there exists a Haar measure on $G$, that is, a
regular invariant Borel measure that is finite for compact sets and positive
for non-empty open sets. This measure, which is unique up to a positive
multiplicative constant, plays a fundamental role in the study of locally
compact groups. Unfortunately, it is known that non-locally compact Polish
groups admit no Haar measure. However, the notion of a Haar nullset has a very
well-behaved generalisation. The following definition was invented by
Christensen \cite{christ}, and later rediscovered by Hunt, Sauer and Yorke
\cite{HSY}. (Actually, Christensen's definition was what we call
generalised Haar null below, but this subtlety will only play a role later.)

\defi{1}{A set $X \subset G$ is called \emph{Haar null} if there exists a Borel
  set $B
\supset X$ and a Borel probability measure $\mu$ on $G$ such that $\mu(B+g)=0$
for every $g \in G$.}

Note that the term \emph{shy} is also commonly used for Haar null, and
co-Haar null sets are often called \emph{prevalent}.

Christensen showed that the Haar null  sets form a
$\si$-ideal, and also that in locally compact groups a set is Haar null iff
it is of measure zero with respect to the Haar
measure. During the last two decades Christensen's notion has been very useful in
studying exceptional sets in diverse areas such as analysis, functional
analysis, dynamical systems, geometric measure theory, group theory, and
descriptive set theory.

Therefore it is very important to understand the fundamental properties of
this $\si$-ideal, such as the Fubini properties, ccc-ness, and all other
similarities and differences between the locally compact and the general case. 

One such example is the following very natural question, which was Problem 1 in
Mycielski's celebrated paper \cite{myc} more than 20 years ago, and was also discussed e.g.
in \cite{dough}, \cite{Banakh} and \cite{sol2}.

\que{3}{\textbf{\textup{[J. Mycielski]}} Let $G$ be a Polish group. Can every Haar null set be covered by a
$G_\delta$ Haar null set?}

It is easy to see using the regularity of Haar measure that the answer is in
the affirmative if $G$ is locally compact.

The first main goal of the present paper is to answer this question.
 
\thm{Main}{If $G$ is a non-locally compact abelian Polish group then there exists
  a (Borel) Haar null set $B \su G$ that cannot be covered by a $G_\delta$ Haar
  null set.}

Actually, the proof will immediately yield that
$G_\delta$ can be replaced by any other class of the Borel hierarchy.
As usual, $\mathbf{\Pi}^0_\xi$ stands for the $\xi$'th multiplicative class
of the Borel hierarchy.

\thm{Maingen}{If $G$ is a non-locally compact abelian Polish group and $1 \le
  \xi < \omega_1$ then there exists a (Borel) Haar null set $B \su G$ that cannot
  be covered by a $\mathbf{\Pi}^0_\xi$ Haar null set.}

It was pointed out to us by Sz. G\l \c ab, see e.g. \cite[Proposition 5.2]{borod}
that an easy but very surprising consequence of this theorem is the following. For the definition of the
additivity of an ideal see e.g. \cite{BJ}.

\cor{ccc}{If $G$ is a non-locally compact abelian Polish group then the
additivity of the $\si$-ideal of Haar null sets is $\om_1$.}

In order to be able to formulate the next question we need to introduce a
slightly modified notion of Haar nullness.
Numerous authors actually use the following weaker definition, in which $B$ is
only required to be universally measurable. (A set is called
\emph{universally measurable} if it is measurable with respect to every Borel
probability measure. Borel measures are identified with their completions.) 

\defi{22}{A set $X \subset G$ is called \emph{generalised Haar null} if there
  exists a universally measurable set $B
\supset X$ and a Borel probability measure $\mu$ on $G$ such that $\mu(B+g)=0$
for every $g \in G$.}

In almost all applications $X$ is actually Borel, so it does not matter which
of the above two definitions we use. Still, it is of substantial theoretical
importance to understand the relation between the two definitions. 
The next question is from Fremlin's problem list \cite{frem}.

\que{3q:frem}{\textbf{\textup{[D. H. Fremlin, Problem GP]}} Is every generalised
Haar null set Haar null? In other words, can every generalised Haar null set be
covered by a Borel Haar null set?}

Dougherty \cite[p.86]{dough} showed that under the Continuum Hypothesis or
Martin's Axiom the answer is in the negative in every non-locally compact Polish group with a two-sided
invariant metric. Later Banakh \cite{Banakh} proved the same under slightly
different set-theoretical assumptions. Dougherty uses transfinite induction,
and Banakh's proof is basically an existence proof using that the so called
cofinality (see e.g. \cite{BJ} for the definition) of the $\si$-ideal of generalised Haar null sets is greater than the continuum
in some models, hence these examples are clearly very far from being Borel.

\medskip

The second main goal of the paper is to answer Fremlin's problem in $ZFC$.

\medskip

Recall that a set is \emph{analytic} if it is the continuous image of
a Borel set, and \emph{coanalytic} if its complement is analytic. Analytic
and coanalytic sets are known to be universally measurable. Since Solecki
\cite{sol2} proved that every analytic generalised Haar null
set is contained in a Borel Haar null set, the following result is optimal.

\thm{Coanal}{Not every generalised Haar null set is Haar null. More precisely, 
if $G$ is a non-locally compact abelian Polish group then there
exists a coanalytic generalised Haar null set $P \su G$ that cannot be covered
by a Borel Haar null set.}

\rem{cptsupp}{We close this section by remarking that in both versions of the
above definition certain authors  actually require
that the measure $\mu$, which we will often refer to as a \emph{witness measure},
has compact support. This is quite important if the underlying group is 
non-separable. However, in our case this would make no difference, since in a
Polish space for every Borel probability measure there exists a compact set
with positive measure \cite[17.11]{cdst}, and then restricting the measure to
this set and normalising yields a witness with a compact support. Therefore we
may suppose throughout the proofs that our witness measures have compact support.}



For more results concerning fundamental properties and applications of Haar null sets
in non-locally compact groups see e.g.~\cite{aubry}, \cite{elekbalk}, \cite{darji},
\cite{dod1}, \cite{dod2}, \cite{elekstep}, \cite{hol}, \cite{mat}, \cite{sol3}, \cite{zaj1}.

\section{Notation and basic facts}
The following notions and facts can all be found in \cite{cdst}. 

Let $\mathcal{F}(G)$ denote the family of closed subsets of $G$ equipped with
the so called Effros Borel structure. Let $\mathcal{K}(G)$ be the family of
compact subsets of $G$ equipped with the Hausdorff metric. Then
$\mathcal{K}(G)$ is a Borel subset of $\mathcal{F}(G)$ and the inherited Borel
structure on $\mathcal{K}(G)$ coincides with the one given by the Hausdorff
metric.

Let us denote
by $\mathcal{P}(G)$ the set of Borel probability measures on $G$, where by
Borel probability measure we mean the completion of a probability measure
defined on the Borel sets. These
measures form a Polish space equipped with the weak*-topology.  For $\mu \in
\mathcal{P}(G)$ we denote by $\supp(\mu)$ the support of $\mu$, i.e. the
minimal closed subset $F$ of $G$ so that $\mu(F)=1$. 
Let $\mathcal{P}_c(G) = \{\mu \in \mathcal{P}(G) : \supp (\mu) \textrm{ is
  compact}\}$.

$\mathbf{\Pi}^0_\xi$ stands for the $\xi$'th multiplicative level of
the Borel hierarchy, $\mathbf{\Delta}^1_1$, $\mathbf{\Sigma^1_1}$ and
$\mathbf{\Pi}^1_1$ denote the classes of Borel, analytic and coanalytic sets,
respectively. For a Polish space $X$, $\mathbf{\Pi}^0_\xi(X)$,
$\mathbf{\Delta}^1_1(X)$ etc.~denote the collections of subsets of $X$ in the
appropriate classes.  Symbols $\Gamma$ and $\Lambda$ will denote one of the
above mentioned classes, and $\check{\Lambda}=\{A^c:A \in \Lambda\}$.

For a set $H \su X \times Y$ we define its $x$-section as $H_x = \{y \in Y :
(x,y) \in H \}$, and similarly if $H \su X \times Y \times Z$ then $H_{x,y} =
\{z \in Z : (x,y,z) \in H \}$, etc. For a function $f \colon X \times Y \to Z$
the $x$-section is the function $f_x \colon Y \to Z$ defined by $f_x(y) =
f(x,y)$. We will sometimes also write $f_x = f(x, \cdot)$.

For $A,B \subset G$ let $d(A,B) = \inf\{d(a,b):a\in A, b\in B\}$ and $A+B
= \{a+b:a \in A, b\in B\}$. Let us denote by $B(g,r)$ and $\bar{B}(g,r)$ the
open and closed ball centered at $g$ of radius $r$.

\section{The proofs}

\subsection{A function with a surprisingly thick graph}

Throughout the proofs, let $\Gamma=\mathbf{\Delta}^1_1$ and
$\Lambda=\mathbf{\Pi}^0_\xi$ for some $1 \leq \xi<\omega_1$, or let
$\Gamma=\mathbf{\Pi}^1_1$ and $\Lambda=\mathbf{\Delta}^1_1$.

The following result will be the starting point of our constructions. For a
fixed measure $\mu$ statement \ref{Bor2}.~below describes the following strange
phenomenon: There exists a Borel graph of a function in a product space such
that every $G_\delta$ cover of the graph has a vertical section of positive
measure.

\thm{Bor}{Let $\Gamma=\mathbf{\Delta}^1_1$ and $\Lambda=\mathbf{\Pi}^0_\xi$
for some $1 \leq \xi<\omega_1$, or let $\Gamma=\mathbf{\Pi}^1_1$ and
$\Lambda=\mathbf{\Delta}^1_1$. Then there exists a partial function
$f:\mathcal{P}_c(G) \times 2^\om \to G$ with $\graph(f) \in \Gamma$ satisfying the
following properties: $\forall \mu \in \mathcal{P}_c(G)$
\begin{enumerate}
 \item\label{Bor1} $(\forall x \in 2^\om)\left[ (\mu, x) \in \dom(f) \Rightarrow f(\mu,x) \in \supp (\mu)\right]$,
 \item\label{Bor2} $(\forall S \in \Lambda(2^\om \times G))\left[(\graph(f_\mu)
     \subset S \Rightarrow (\exists x \in 2^\om)(\mu(S_x)>0 )\right].$
\end{enumerate} }

Before the proof we need several technical lemmas.

\lem{Pc}{$\mathcal{P}_c(G)$ is a Borel subset of $\mathcal{P}(G)$.}

\proof
The map $\mu \mapsto \supp(\mu)$ between $\mathcal{P}(G)$ and $\mathcal{F}(G)$ is
Borel (see \cite[17.38]{cdst}) and $\mathcal{P}_c(G)$ is the preimage of
$\mathcal{K}(G)$ under this map.
\qed

\lem{tek}{Let $X$ be a Polish space and $C \subset \mathcal{P}_c(G) \times X
\times G$ with $C \in \Gamma$. Then $\{(\mu,x):\mu(C_{\mu,x})>0\} \in \Gamma$.}

\proof
Let first $\Gamma = \mathbf{\Delta}^1_1.$
If $Y$ is a Borel space and $C \subset Y \times G $ is a Borel set then the
map $\varphi \colon Y \times \mathcal{P}_c(G) \to [0,1]$ defined by
$\varphi(y,\mu) = \mu(C_y)$ is Borel (\cite[17.25]{cdst}). Using this
for $Y=\mathcal{P}_c(G) \times X$ we obtain that the map $\psi \colon
\mathcal{P}_c(G) \times X \to [0,1]$ given by $\psi(\mu, x) = \varphi((\mu,
x), \mu) = \mu(C_{\mu,x})$  is also Borel. Then $\{(\mu,x):\mu(C_{\mu,x})>0\} =
\psi^{-1}((0,1])$, hence Borel.

For $\Gamma = \mathbf{\Pi}^1_1$ this is simply a special case of \cite[36.23]{cdst}.
\qed



\lem{tek1}{The set $\{(\mu,g):g \in \supp(\mu)\} \subset \mathcal{P}_c(G) \times G$ is Borel.}

\proof
As mentioned above, the map $\mu \mapsto \supp(\mu)$ is Borel between
$\mathcal{P}(G)$ and $\mathcal{F}(G)$, hence its restriction to
$\mathcal{P}_c(G)$ is also Borel.
 
Let $E=\{(K,g):K \in \mathcal{K}(G),\  g \in K\}$, which clearly is a closed
subset of $\mathcal{K}(G) \times G$. If we denote by $\Psi:\mathcal{P}_c(G)
\times G \to \mathcal{K}(G) \times G$ the Borel map defined by $(\mu,g) \mapsto
(\supp(\mu),g)$ then we obtain that $\{(\mu,g):g \in
\supp(\mu)\}=\Psi^{-1}(E)$ is Borel.
\qed

Let us now prove Theorem \ref{Bor}.

\bigskip

\proof
Let $U \in \Gamma(2^\omega \times 2^\om \times G)$ be
 universal for the $\check{\Lambda}$ subsets of $2^\om \times
G$, that is, for every $A \in \check{\Lambda}(2^\om \times
G)$ there exists an $x \in 2^\omega$ such that $U_x = A$ (for the existence of
such a set see \cite[22.3, 26.1]{cdst}). Notice that $\check{\Lambda} \subset
\Gamma$. Let
\[
U'=\mathcal{P}_c(G) \times U.
\]
Define
\[
U''=\{(\mu,x,g) \in \mathcal{P}_c(G) \times 2^\omega \times G: (\mu,x,x,g) \in
U' \text{ and }\mu(U'_{\mu,x,x})>0\},
\]
then $U'' \in \Gamma$ using that the map $(\mu, x, g) \mapsto (\mu, x, x, g)$ is continuous and by Lemma
\ref{tek}.
Let
\[
U'''=\{(\mu,x,g) \in U'': g \in \supp(\mu)\},
\] 
then $U''' \in \Gamma$ by Lemma \ref{tek1}. Clearly,
\[U'''_{\mu,x}=
\begin{cases}
  U'_{\mu,x,x} \cap \supp(\mu) & \text{if } \mu(U'_{\mu,x,x})>0, \\
  \emptyset & \text{otherwise.} 
\end{cases}
\]

Since for all $(\mu,x)$ the section $U'''_{\mu,x}$ is either empty or has positive $\mu$ measure, by the 'large
section uniformisation theorem' \cite[18.6]{cdst} and the coanalytic uniformisation theorem \cite[36.14]{cdst} there exists a 
partial function $f$ with $\graph(f) \in \Gamma$ such that $\dom(f) = \{(\mu,x) \in \mathcal{P}_c(G) \times 2^\omega : \mu(U'_{\mu,x,x})>0\}$ and $\graph(f) \subset U'''$.

We claim that this $f$ has all the required properties. 

First, by the definition of $U'''$, clearly $f(\mu,x) \in
\supp(\mu)$ holds whenever $(\mu,x) \in \dom(f)$, hence Property \ref{Bor1}.~of Theorem \ref{Bor} holds.

Let us now prove Property \ref{Bor2}. Assume towards a contradiction that there exists $\mu \in \mathcal{P}_c(G)$ and $S \in \Lambda(2^\om \times G)$ such that $\graph(f_\mu) \subset S$ and $\mu(S_x)=0$ for every $x \in 2^\omega$. Define $B=(2^\om \times G) \setminus S$. By the universality of $U$
there exists $x \in 2^\omega$ such that $U_x = U'_{\mu,x} = B$.
Now, for every $y \in 2^\om$ the section $B_{y}$ is of positive (actually full) $\mu$ measure,
in particular $\mu(U'_{\mu,x,x})>0$, and therefore $(\mu,x) \in \dom(f)$ and
\[
f(\mu,x) \in U'''_{\mu,x}\subset U''_{\mu,x}=U'_{\mu,x,x}=B_x.
\]
However, $f(\mu,x) \in S_x=G \setminus B_x$, a contradiction.
\qed

\subsection{Translating the compact sets apart}

This section heavily builds on ideas of Solecki \cite{sol1}, \cite{sol2}. The
main point is that if $G$ is non-locally compact then one can apply a
translation (chosen in a Borel way) to every compact subset of $G$ so that the
resulting translates are disjoint. (For technical reasons we will need to 
consider continuum many copies of each compact set and also to `blow them up' 
by a fixed compact set $C$.)

\prop{eltol}{Let $C \in \mathcal{K}(G)$ be fixed. Then there exists a Borel map $t:
\mathcal{K}(G) \times 2^\om \times 2^\om \to G$ so that
\begin{enumerate}
\item\label{eltol1} if $(K,x,y) \not = (K',x',y')$ are elements of $\mathcal{K}(G) \times 2^\om \times 2^\om$ then
\[
(K - C + t(K,x,y)) \cap (K' - C' + t(K',x',y'))= \emptyset
\]
\item\label{eltol2} for every $K \in \mathcal{K}(G)$ and $y \in 2^\om$  the
  map $t(K,\cdot, y)$ is continuous.
\end{enumerate}}

\proof
We use Solecki's arguments \cite{sol1}, \cite{sol2}, which he used for
different purposes, with some modifications. However, for the sake of
completeness, we repeat large parts of his proofs.

Fix an increasing sequence of finite sets $Q_k \su G$ with $0 \in Q_0$ such
that $\cup_{k \in \om} Q_k$ is dense in $G$.

\lem{eltol?}{For every $\epsilon>0$ there exists $\delta>0$ and a sequence
$\{g_k\}_{k \in \om} \su B(0,\epsilon)$ such that for every distinct $k, k'
\in \omega $
\[ d( Q_k+g_k, Q_{k'} +g_{k'}) \ge \delta.
\]}

\proof
Since $G$ is not locally compact, there
exists $\delta>0$ and a countably infinite set $S \subset B(0,\epsilon)$ such
$d(s,s') \ge 2\de$ for every distinct $s, s' \in S$. 

Now we define $g_k$ inductively as follows. Suppose that we are done for
$i<k$. If for every $s \in S$ there are $a \in Q_k$, $i < k$ and $b \in Q_i$ with
$d(a+s, b+g_i) < \delta$ then there is a pair $s,s'$ of distinct members
of $S$ with the same $a, i$ and $b$. But then   
\[
d(s, s') = d(a+s, a+s' ) \le d(a+s, b+g_i ) + d(b + g_i, a+s') < 2 \delta,
\]
a contradiction. Hence we can let $g_k = s$ for an appropriate $s \in S$.
\qed

It is easy to see that using the previous lemma repeatedly we can inductively
fix $\epsilon_n$, $\de_n <\eps_n$ and sequences $\{ g^n_k\}_{k \in \om}$ such
that for every $n \in \om$
\begin{itemize}
\item
$\{ g^n_k\}_{k \in \om} \su B(0, \eps_n)$,
\item
$d(Q_k+g^n_k, Q_{k'}+g^n_{k'} ) \ge 2 \delta_n$ for every distinct $k, k'
\in \om$,
\item
$\sum_{m>n} \epsilon_{m} < \frac{\delta_n}{3}.$
\end{itemize}
 
Note that the second property implies that for every $n \in \om$ the function
$k \mapsto g^n_k$ is injective. Note also that $\eps_n \to 0$ and hence $\de_n
\to 0$, moreover, $\sum \de_n$ is also convergent.

Let us also fix a Borel injection $c:\mathcal{K}(G) \times 2^\om \times 2^\om
\to \om^\om$ such that for each $K$ and $y$ the map $c(K,\cdot,y)$ is
continuous. (E.g. fix a Borel injection $c_1 : \mathcal{K}(G) \to
2^\omega$ and continuous injection $c_2:2^\om \times 2^\om \times 2^\om \to
\om^\om$ and let $c(K,x,y)=c_2(c_1(K),x,y)$.)

Our goal now is to define $t(K,x,y)$, so let us fix a triple
$(K,x,y)$. First we define a sequence $\{h_n = h_n(K,x,y)\}_{n \in \om}$ with
$h_n \in \{g^n_k\}_{k \in \om}$ as follows.
Suppose that we are given $h_i$ for $i<n$. By
the density of $\cup_k Q_k$ we have $G=\cup_k (Q_{k}+B(0,\delta_n/2))$.
Since $K-C$ is compact, there exists a minimal index $k_n(K,x,y)$ so that
\[
K - C + \sum_{i<n} h_i \subset Q_{k_n(K,x,y)}+B(0,\delta_n/2).
\]
Fix an injective map $\phi : \om \times \om \to \om$ with $\phi(i,j) \ge i$ for every $i \in \om$ and let
\begin{equation}
\label{e:hn}
h_n=g^n_{\phi(k_n(K,x,y), c(K,x,y)(n))}
\end{equation}
and 
\begin{equation}
\label{e:t}
t(K,x,y)=\sum_{n \in \om} h_n.
\end{equation} 
We claim that this function has the required properties.

First, it is well defined, that is, the sum is convergent since $h_n \in
B(0,\epsilon_n)$, and hence for all $n \in \om$
\begin{equation}
\label{e:hn2}
\sum_{m>n} h_m \in \bar{B}(0,\delta_n/3).
\end{equation}

In order to prove \ref{eltol1}.~of the Proposition, let us now fix $(K,x,y)
\not = (K',x',y')$. Then
there exists an $n \in \omega$ such that $c(K,x,y)(n) \not = c(K',x',y')(n)$. 
By the injectivity of $\phi$ and of the sequence $k \mapsto g^{n}_k$ and
also by \eqref{e:hn} we obtain that
$h_{n}(K,x,y) \neq h_{n}(K',x',y')$. Denote by $h_i$ and $h'_i$ the elements
$h_i(K,x,y)$ and $h_{i}(K',x',y')$, respectively. Set
\[
k = \phi(k_{n}(K,x,y),c(K,x,y)(n)) \textrm{ and } k' = \phi(k_n(K',x',y'),c(K',x',y')(n)).
\]

The condition $\phi(i,j) \ge i$
implies $k \ge k_n(K,x,y)$, hence $Q_k \supset Q_{k_n(K,x,y)}$ and  similarly
$k' \ge k_n(K',x',y')$, so $Q_{k'} \supset Q_{k_n(K',x',y')}$. Therefore, by
the definition of $k_n$,
\[
K - C + \sum_{i<n} h_i \in Q_k+B(0,\delta_n/2) \textrm{ and } K' - C + \sum_{i<n} h'_i   \in Q_{k'}+B(0,\delta_n/2),
\]
hence
\[
K - C + \sum_{i\le n} h_i \in Q_k + h_n + B(0,\delta_n/2) \textrm{ and } K' - C + \sum_{i\le n} h'_i   \in Q_{k'} + h'_n + B(0,\delta_n/2).
\]
Thus, using the triangle inequality and the second property of the $g^n_k$ we obtain
\[
d(K-C+\sum_{i\leq n} h_i,K'-C+\sum_{i \leq n} h'_i)\geq
d(Q_k+h_n,Q_{k'}+h'_n)-2\cdot \frac{\delta_n}{2}=
\]
\[
= d(Q_k+g^n_k,Q_{k'}+g^n_{k'}) - \delta_n \geq 2\delta_n-\delta_n = \delta_n.
\]
From this, using \eqref{e:hn2}, we obtain $d(K-C+t(K,x,y),K'-C+t(K',x',y'))\ge  \de_n - 2 \frac{\de_n}{3} = \frac{\delta_n}{3} > 0$, which proves \ref{eltol1}.

What remains to show is that $t$ is a Borel map and for every $K$ and $y$ the
map $t(K, \cdot, y)$ is continuous. But \eqref{e:hn2} shows that the series defining $t$ in \eqref{e:t} is uniformly convergent, so the next lemma finishes the proof.

\lem{44}{For every $n \in \om$ the map $h_n$ is Borel and for every $K$ and $y$ the
map $h_n(K, \cdot, y)$ is continuous.}

\proof
We will actually prove more by induction on $n$. 
Define $f_n \colon \mathcal{K}(G) \times 2^\om \times 2^\om \to \mathcal{K}(G)$ by 
\begin{equation}
\label{e:fn}
f_n(K,x,y) = K - C + \sum_{i<n} h_i (K,x,y).
\end{equation}
We claim that the maps $f_n$, $k_n$ and $h_n$ are Borel and for every $K$ and
$y$ the maps $f_n(K, \cdot, y)$, $k_n(K, \cdot, y)$ and $h_n(K, \cdot, y)$ are 
locally constant.

Note that if a function takes its values from a discrete set than locally
constant is equivalent to continuous.

First we prove that the maps are Borel. Suppose that we are done for $i <
n$. Let us check that $f_n$ is Borel. Put $\eta : (K,x,y) \mapsto (K,
\sum_{i<n} h_i (K,x,y))$ and $\psi : (K,g) \mapsto K - C + g$, then $f_n =
\psi \circ \eta$. Moreover, $\eta$ is Borel by induction, and $\psi$ is easily
seen to be continuous, hence $f_n$ is Borel.

Next we show that $k_n$ is Borel. Since $\ran(k_n) \su \om$, 
we need to check that for every fixed $m \in \om$ the set $B = \{(K, x, y)
\colon k_n(K,x,y) = m \}$ is Borel. By the definition of $k_n(K,x,y)$, clearly
\[
B = \{(K,x,y) \colon f_n(K,x,y) \subset U \textrm{ and } f_n(K,x,y) \not\subset V\},
\]
where $U = Q_m + B(0,\delta_n/2)$ and $V = Q_{m-1} + B(0,\delta_n/2)$ are fixed open sets. 

Set $\mathcal{U}_W = \{L \in \mathcal{K}(G): L \subset W\}$, which is open in $\mathcal{K}(G)$ for every open set $W \su G$. Then clearly
\[
B = f_n^{-1} (\mathcal{U}_U) \sm f_n^{-1}(\mathcal{U}_V),
\]
hence Borel.

Since the functions $k \mapsto g^n_k$ and $\phi$ defined on countable sets are
clearly Borel, the Borelness of $k_n$ and $c$ imply by \eqref{e:hn}
that $h_n$ is also Borel.

In order to prove that $f_n$, $k_n$ and $h_n$ are locally constant in the
second variable, fix $K$ and $y$ and suppose that we are done for $i<n$. Then
\eqref{e:fn} shows that $f_n$ is locally constant in the second variable
by induction. This easily implies using the definition of $k_n$ that $k_n$ is
also is locally constant in the second variable. But from this, and from the
fact that $c(K, \cdot, y)(n) \colon 2^\om \to \om$ is continuous, hence
locally constant, it is also
clear using \eqref{e:hn} that $h_n$ is also locally constant in the second
variable, which finishes the proof of the Lemma.
\qed

Therefore the proof of the Proposition is also complete.
\qed

\section{Putting the ingredients together}

Now we are ready to prove our main results, which are summarised in the following theorem.

\thm{Technical}{Let $\Gamma=\mathbf{\Delta}^1_1$ and
$\Lambda=\mathbf{\Pi}^0_\xi$ for some $1 \leq \xi<\omega_1$, or let
$\Gamma=\mathbf{\Pi}^1_1$ and $\Lambda=\mathbf{\Delta}^1_1$.  If $G$ is a
non-locally compact abelian Polish group then there exists a (generalised, in
the case of $\Gamma=\mathbf{\Pi}^1_1$) Haar null set $E \in \Gamma(G)$ that is
not contained in any Haar null set $H \in \Lambda(G)$.}

\proof Let $f$ be given by Theorem \ref{Bor}.

Denote the Borel map $\mu \mapsto \supp(\mu)$ by
$\supp:\mathcal{P}_c(G) \to \mathcal{K}(G)$. 
Let us also fix a Borel bijection $c:\mathcal{P}_c(G) \to 2^\omega$ (which we
think of as a coding map) and a continuous probability measure
$\nu$ on $G$ with compact support $C$ containing $0$ (compactly supported continuous measures exist on every Polish space without isolated points, since such spaces contain copies of $2^\omega$). Let $t:\mathcal{K}(G) \times 2^\om \times
2^\om \to G$ be the map from Proposition \ref{eltol} with the $C$ fixed above, and define the map $\Psi \colon \mathcal{P}_c(G) \times 2^\om \times
G \to G$ by
\begin{equation}
\label{e:psi}
\Psi(\mu,x,g)=g + t(\supp(\mu),x,c(\mu)).
\end{equation}
Finally, define $E=\Psi(\graph(f))$.

\clm{c1}{$E \in \Gamma$.}

\proof
$\Psi$ is clearly a Borel map.  We claim that it is injective on
$D = \{(\mu,x,g):\mu \in \mathcal{P}_c(G), g\in \supp(\mu)\}$, which is Borel by Lemma \ref{Pc} and Lemma \ref{tek1}. Let $(\mu, x, g) \neq (\mu', x', g')$ be elements of $D$, we need to check that $\Psi$ takes distinct values on them. The case $(\mu, x) = (\mu', x')$ is obvious, while the case $(\mu, x) \neq (\mu', x')$ follows from Property \ref{eltol1}.~in Proposition \ref{eltol}, since $\Psi(\mu,x,g) \in \supp(\mu) - C + t(\supp(\mu),x,c(\mu))$ (recall that $g \in \supp(\mu)$ and $0 \in C$).
Therefore $\Psi$ is a Borel isomorphism on $D$. By $\graph(f) \su D$ this
implies that $E=\Psi(\graph(f))$ is in $\Gamma$ (for
$\Gamma=\mathbf{\Delta^1_1}$ see \cite[15.4]{cdst}, for
$\Gamma=\mathbf{\Pi}^1_1$ notice that by \cite[25.A]{cdst} a Borel isomorphism
takes analytic sets to analytic sets, hence coanalytic sets to coanalytic sets).
\qed

\clm{c2}{$E$ is Haar null (generalised Haar null in the case of $\Gamma = \mathbf{\Pi}^1_1$).}

\proof
We prove that $\nu$ is witnessing this fact. Actually, we prove more: $|C \cap (E+g)|\leq 1$ for every $g \in G$, or equivalently
$|(C+g) \cap E|\leq 1$ for every $g \in G$. So let us fix $g \in G$.
\[
E = \Psi(\graph(f)) = \{ \Psi(\mu, x, f(\mu, x)) : (\mu, x) \in \dom(f) \} = 
\]
\[
\{f(\mu, x) + t(\supp(\mu), x, c(\mu)) : (\mu, x) \in \dom(f) \},
\]
hence the elements of $E$ are of the form $g^{\mu,x} = f(\mu, x) +
t(\supp(\mu), x, c(\mu))$. This element $g^{\mu,x}$ is clearly in $A^{\mu,x} =
\supp(\mu) + t(\supp(\mu), x, c(\mu))$ by Property \ref{Bor1}.~of
Theorem \ref{Bor}, and the sets $A^{\mu,x}$ form a pairwise disjoint
family as $(\mu, x)$ ranges over $\dom(f)$, by Property \ref{eltol2}.~of
Proposition \ref{eltol}. Hence it suffices to show that $C+g$ can intersect at
most one $A^{\mu,x}$. But it can actually intersect at most one set of the
form $K + t(K, x, y)$, since otherwise $g$ would be in the intersection of two
distinct sets of the form $K - C + t(K, x, y)$, contradicting Property
\ref{eltol2}.~of Proposition \ref{eltol}.
\qed

\clm{c3}{There is no Haar null set $H \in \Lambda$ containing $E$.}

Suppose that $H \in \Lambda$ is such a set. Then by Remark \ref{cptsupp} there exists a probability measure $\mu$ 
with compact support witnessing this fact. The section map $\Psi_{\mu} = \Psi(\mu,\cdot,\cdot)$ is continuous by \eqref{e:psi} and Property \ref{eltol2}.~of Proposition \ref{eltol}. Now let $S=\Psi^{-1}_{\mu}(H)$, then $S \in \Lambda(2^\om \times G)$.

It is easy to check that $\graph(f_\mu) \subset S$, and therefore, using
Theorem \ref{Bor}, there exists $x \in 2^\om$ such that $\mu(S_{x})>0$. By the
definition of $S$ we have that $\Psi(\mu,x,S_x) \subset \Psi_\mu(S) \subset
H$. But $\Psi(\mu,x,\cdot):G \to G$ is a translation, so a translate of $H$
contains $S_x$, which is of positive $\mu$ measure, contradicting that $H$ is
Haar null with witness $\mu$.
\qed

This concludes the proof.
\qed

\section{Questions}

\que{q1}{Let $G$ be a non-locally compact abelian Polish group. Does there exist an $F_\sigma$ Haar null set that cannot be covered by a $G_\delta$ Haar null set?}

Interestingly, our proof does not give any information about the Borel class of our example.

\que{q2}{What is the least complexity of such a set? And in general, what is the least complexity of a Haar null set that cannot be covered by a $\mathbf{\Pi}^0_\xi$ Haar null set?}

\rem{rrr}{We remark that it is not hard to show that in abelian Polish groups every $\si$-compact Haar null set can be covered by a $G_\delta$ Haar null set.}

\que{q3}{Do the results of the paper hold in all (not necessarily abelian) non-locally compact Polish groups?}

\que{q4}{Does there exist a Polish group with a countable subset that cannot be covered by a $G_\delta$ Haar null set?}

In view of the above remark, the group in the last question cannot be abelian. Of course, it also cannot be locally compact. How about e.g.~an arbitrary countable dense subset of $Homeo[0,1]$? This is actually closely related to the following question, popularised by U. B. Darji, and considered e.g. in \cite{cohen}.

\que{q5}{Can every uncountable Polish group be written as a union of a Haar null set and a meager set?}

The answer is affirmative e.g. for abelian groups or for groups with a two-sided invariant metric.

\medskip

The so called cardinal invariants convey a lot of information about the set-theoretical properties of a $\si$-ideal, see e.g. \cite{BJ}. Banakh examined this problem in detail in \cite{Banakh} for the $\si$-ideal of \emph{generalised} Haar null sets.

\que{q6}{What can we say about the cardinal invariants of the $\si$-ideal of Haar null sets? How about e.g. if $G = \ZZ^\omega$?}

Surprisingly, the invariants may differ for Haar null and generalised Haar null sets. First, in contrast with Corollary \ref{ccc}, \cite[Thm. 3]{Banakh} shows that the additivity of the generalised Haar null sets in $\ZZ^\omega$ equals the additivity of the Lebesgue null sets. Second, \cite[Thm. 3]{Banakh} also shows that the cofinality of the generalised Haar null sets in $\ZZ^\omega$ may exceed the continuum, whereas for Haar null sets it is clearly at most continuum.

\medskip

In separable Banach spaces there is a well-known alternative notion of nullness. For the equivalent definitions of Aronszajn null, cube null and Gaussian null sets see \cite{mari}.

\que{q7}{Suppose that $G$ is a separable Banach space. Which results of the paper remain valid when Haar null is replaced by Aronszajn null?}

\bigskip

M\'arton Elekes

Alfr\'ed R\'enyi Institute of Mathematics

Hungarian Academy of Sciences

P.O. Box 127, H-1364 Budapest, Hungary

elekes.marton@renyi.mta.hu

www.renyi.hu/ $\tilde{}$ emarci

and

E\"otv\"os Lor\'and University

Department of Analysis

P\'az\-m\'any P. s. 1/c, H-1117, Budapest, Hungary

\bigskip

Zolt\'an Vidny\'anszky

E\"otv\"os Lor\'and University

Department of Analysis

P\'az\-m\'any P. s. 1/c, H-1117, Budapest, Hungary

www.cs.elte.hu/ $\tilde{}$ vidnyanz


\begin{thebibliography}{abc}


\bibitem{aubry} J. Aubry, F. Bastin, S. Dispa, Prevalence of multifractal functions in $S^\nu$ spaces, \textit{J. Fourier Anal. Appl.} \textbf{13} (2007), no. 2, 175--185.

\bibitem{elekbalk} R. Balka, U. B. Darji, M. Elekes, Bruckner-Garg-type results with respect to Haar null sets in C[0,1], submitted.

\bibitem{Banakh} T. Banakh, Cardinal characteristics of the ideal of Haar null sets, \textit{Comment. Math. Univ. Carolinae} \textbf{45} (2004),  no. 1, 119--137.

\bibitem{BJ} T. Bartoszy\'nski and H. Judah, \textsl{Set theory. On the structure
of the real line.} A K Peters, Ltd., Wellesley, MA, 1995.

\bibitem{kechbeck} H. Becker, A. S. Kechris, \textit{The descriptive set theory of Polish group actions}. London Mathematical Society Lecture Note Series, 232. Cambridge University Press, Cambridge, 1996.

\bibitem{borod} P. Borodulin-Nadzieja, Sz. G\l \c ab, Ideals with bases of unbounded Borel complexity, \textit{MLQ Math. Log. Q.} \textbf{57} (2011), no. 6, 582--590.

\bibitem{christ} J. P. R. Christensen, On sets of Haar measure zero in abelian Polish groups, \textit{Israel J. Math.} \textbf{13} (1972), 255--260. 

\bibitem{cohen} M. P. Cohen, R. R. Kallman, Openly Haar null sets and conjugacy in Polish groups, preprint.

\bibitem{mari} M. Cs\"ornyei, Aronszajn null and Gaussian null sets coincide, \textit{Israel J. Math.} \textbf{111} (1999), no. 1, 191--201.

\bibitem{darji} B. U. Darji, On Haar meager sets, \textit{Topology Appl.} \textbf{160} (2013), no. 18, 2396--2400.

\bibitem{dod1} P. Dodos, On certain regularity properties of Haar-null sets, \textit{Fund. Math.} \textbf{181} (2004), no. 2, 97--109. 

\bibitem{dod2} P. Dodos, The Steinhaus property and Haar-null sets, \textit{Bull. Lond. Math. Soc.} \textbf{41} (2009), 377--384.

\bibitem{dough} R. Dougherty, Examples of non-shy sets, \textit{Fund. Math.} \textbf{144} (1994), 73-88. 

\bibitem{elekstep} M. Elekes, J. Stepr\=ans, Haar null sets and the consistent reflection of non-meagreness, \textit{Canad. J. Math.} \textbf{66} (2014), 303--322.

\bibitem{frem} D. H. Fremlin, Problems,\\ http://www.essex.ac.uk/maths/people/fremlin/problems.pdf.


\bibitem{hol} P. Holick\'y, L. Zaj\'i\v cek, Nondifferentiable functions, Haar null sets and Wiener measure, \textit{Acta Univ. Carolin. Math. Phys.} \textbf{41} (2000), no. 2, 7--11.

\bibitem{HSY} B. Hunt, T. Sauer, J. Yorke, Prevalence: a translation-invariant ``almost
every'' on infinite-dimensional spaces, \textit{Bull. Amer. Math. Soc.} \textbf{27} (1992), 217--238.

\bibitem{cdst} A. S. Kechris, \textit{Classical Descriptive Set
Theory}, Graduate Texts in Mathematics 156, Springer-Verlag, New York, 1995.

\bibitem{mat} E. Matou\v skov\'a, The Banach-Saks property and Haar null sets, \textit{Comment. Math. Univ. Carolin.} \textbf{39} (1998), no. 1, 71--80.

\bibitem{myc} J. Mycielski, Some unsolved problems on the prevalence of ergodicity, instability, and algebraic independence, \textit{Ulam Quarterly} \textbf{1} (1992), 30--37.

\bibitem{sol1} S. Solecki, Haar null and non-dominating sets, \textit{Fund. Math.} \textbf{170} (2001), 197--217.

\bibitem{sol2} S. Solecki, On Haar null sets, \textit{Fund. Math.} \textbf{149} (1996), 205--210.

\bibitem{sol3} S. Solecki, Size of subsets of groups and Haar null sets, \textit{Geom. Funct. Anal.} \textbf{15} (2005), no. 1, 246--273.

\bibitem{zaj1} L. Zaj\'i\v cek, On differentiability properties of typical continuous functions and Haar null sets, \textit{Proc. Amer. Math. Soc.} \textbf{134} (2006), no. 4, 1143--1151.

\end{thebibliography}
\end{document}